\documentstyle{article}

\newcommand{\N}{{\rm I\kern-.5ex N}}
\newcommand{\Z}{{\sf \vrule height 1.55ex depth-1.2ex width.03em\kern-.11em Z \kern-.9ex Z\kern-.11em\vrule height 0.3ex depth0ex width.03em}}
\newcommand{\Q}{{\rm\kern.2ex\vrule height1.55ex depth-.05ex width.03em\kern-.7ex Q}}
\newcommand{\R}{{\rm I\kern-.5ex R}}
\newcommand{\Rvar}{{\rm I\kern-.5ex R}}
\newcommand{\C}{{\rm\kern.3ex\vrule height1.55ex depth-.05ex width.03em\kern-.7ex C}}
\newcommand{\Cvar}{{\, \rm\kern.3ex\vrule height1.1ex depth-.05ex width.03em\kern-.7ex C}}

\newcommand{\spat}{\hspace{4ex}}

\newcommand{\op}{B(H)}
\newcommand{\cop}{B_0(H)}

\newcommand{\ar}{A_r}

\newcommand{\nab}{\nabla}

\newcommand{\cN}{{\cal N}}
\newcommand{\cM}{{\cal M}}

\newcommand{\od}{\odot}
\newcommand{\ot}{\otimes}

\newcommand{\la}{\Lambda}

\newcommand{\om}{\omega}
\newcommand{\io}{\iota}
\newcommand{\vfi}{\varphi}
\newcommand{\vep}{\varepsilon}

\newcommand{\al}{\alpha}

\newcommand{\de}{\Delta}

\newcommand{\th}{\theta}

\newcommand{\si}{\sigma}
\newcommand{\Mfi}{{\cal M}_{\vfi}}
\newcommand{\Nfi}{{\cal N}_{\vfi}}

\newcommand{\lafi}{\la_\vfi}

\newcommand{\pifi}{\pi_\vfi}

\newcommand{\Nfir}{{\cal N}_{\vfi_r}}
\newcommand{\Mfir}{{\cal M}_{\vfi_r}}

\newcommand{\lar}{\la_r}

\newcommand{\text}[1]{\mbox{#1}}

\newcommand{\cst}{\text{C}$^*$}

\newcommand{\qed}{\ \hfill \rule{2mm}{2mm}}
\newenvironment{demo}{\medskip\noindent\bf Proof :\ \  \rm}{\qed\bigskip\par }

\newtheorem{definition}{Definition}[section]
\newtheorem{proposition}[definition]{Proposition}
\newtheorem{lemma}[definition]{Lemma}

\newtheorem{remark}[definition]{Remark}
\newtheorem{theorem}[definition]{Theorem}

\newtheorem{result}[definition]{Result}

\begin{document}
\begin{center}
\LARGE\bf
A natural extension of a left invariant lower semi-continuous weight
\end{center}

\bigskip

\begin{center}
\rm J. Kustermans  \footnote{Research Assistant of the
National Fund for Scientific Research (Belgium)}

Institut for Matematik og Datalogi

Odense Universitet

Campusvej 55

5230 Odense M

Denmark

\bigskip\medskip

\bf April 1997 \rm
\end{center}

\bigskip

\subsection*{Abstract}
In this paper, we describe a natural method to extend left invariant weights on
\cst-algebraic quantum groups. This method is then used to improve the invariance property of a
left invariant weight. We also prove some kind of uniqueness result for left Haar weights
on \cst-algebraic quantum groups arising from algebraic ones.

\section*{Introduction}

Possibly the most important object associated to locally compact groups is the left Haar
measure. So it is not a great surprise that a major role in the \cst-algebraic approach to
quantum groups is played by so called left Haar weights. However, there are some different
ways to expres the left invariance.

We will discuss one of them in this introduction.

\medskip

Therefore, consider a Hilbert space $H$ and a non-degenerate sub-\cst-algebra $B$ of
$B(H)$.

Let $v,w \in H$, then $\om_{v,w}$ denotes the element in $B^*$ such that
$\om_{v,w}(x) = \langle x \, v , w \rangle$ for every $x \in B$.

Consider moreover a comultiplication $\de$ on $B$, this is a non-degenerate
$^*$-homomorphism from $B$ into $M(B \ot B)$ which is coassociative: $(\de \ot \io)\de  =
(\io \ot \de)\de$.

\medskip

Probably the weakest interesting form of left invariance is the following one.

Let $\vfi$ be a densely defined lower semi-continuous weight on $B$, then we call $\vfi$
weakly left invariant if and only if we have for every $x \in B^+$ with $\vfi(x) < \infty$
and every $v \in H$ that $$ \vfi\bigl( (\om_{v,v} \ot \io)\de(x) \bigr)  = \langle v , v
\rangle \, \vfi(x) . $$

It seems to be an interesting question whether it is  reasonable to require that a left
invariant weight satisfies also some kind of converse :

Let $\vfi$ be a densely defined lower semi-continuous weight on $B$, then we call $\vfi$
left invariant if and only if $\vfi$ satisfies  the following two conditions :
\begin{enumerate}
\item $\vfi$ is weakly left invariant
\item Let $x$ be an element in $B^+$. If $\vfi\bigl((\om_{v,v} \ot \io)\de(x)\bigr) <
\infty$ for every $v \in H$, then we have that $\vfi(x) < \infty$.
\end{enumerate}

\medskip

We prove in the first section that under some extra conditions on $(B,\de)$, any weakly
left invariant weight has a smallest extension which is left invariant. We would like to
mention that \cst-algebraic quantum groups according to an upcoming definition of Masuda,
Nakagami and Woronowicz will satisfy these extra conditions.

\medskip

In the last section, we apply the results of the first section to prove two results
concerning left invariant weights on \cst-algebraic quantum groups arising from algebraic
ones (see \cite{Kus} and \cite{VD1}) :
\begin{itemize}
\item We show that the canonical left Haar weight  is automatically left invariant in the sense above.
As a consequence, we get also that the form of left invariance mentioned above is
equivalent with an ever stronger form of left invariance involving the slice
\cst-valued weight $\io \ot \vfi$.
\item We prove that left invariant weights are uniquely determined under some relative invariance condition.
A similar result for quantum E(2) can be found in \cite{Baa}.
\end{itemize}

It seems that the techniques used in this last section can also be applied to the more
general case of \cst-algebraic quantum groups according to Masuda, Nakagami \& Woronowicz.

\bigskip

We should mention that the used terminology of left invariance and weak left invariance is
by no means standard. It was only introduced to be able to talk about left invariance.

\bigskip

We end this introductions with some basic notations and results.

\medskip

For any Hilbert space $H$, the set of bounded operators on $H$ will be denoted by $B(H)$
whereas the set of compact operators on $H$ will be denoted by $B_0(H)$.

\medskip

Consider a \cst-algebra $B$. A norm continuous one-parameter group $\al$ on $B$ is a group
homomorphism from $\R$ into the group of $^*$-automorphisms on $B$ such that the mapping
$\R \rightarrow B : b \mapsto \al_t(b)$ is norm continuous for every $b \in B$.

Using the theory of analytic functions, we can define for every complex number $z$ a
closed linear operator $\al_z$ from within $B$ into $B$. (see e.g. \cite{Zsido} )

\medskip

Let $\vfi$ be a densely defined lower semi-continuous weight on a \cst-algebra $B$. We will use the
following notations:
\begin{itemize}
\item ${\cal M}^+_\varphi = \{\,  b \in B^+ \mid \varphi(b) < \infty  \,\} $
\item ${\cal N}_\varphi = \{\,  b \in B \mid \varphi(b^*b) < \infty \,\} $
\item ${\cal M}_\varphi = \text{span\ } {\cal M}^+_\varphi
= {\cal N}_\varphi^* {\cal N}_\varphi$ .
\end{itemize}

\medskip

A GNS-construction of $\varphi$ is by definition a triple
$(H_\varphi,\pi_\varphi,\Lambda_\varphi)$ such that
\begin{itemize}
\item $H_\varphi$ is a Hilbert space
\item $\Lambda_\varphi$ is a linear map from ${\cal N}_\varphi$ into
      $H_\varphi$ such that
      \begin{enumerate}
      \item  $\Lambda_\varphi({\cal N}_\varphi)$ is dense in $H_\varphi$
      \item   We have for every $a,b \in {\cal N}_\varphi$ that
              $\langle \Lambda_\varphi(a),\Lambda_\varphi(b) \rangle
              = \varphi(b^*a) $
      \end{enumerate}
      Because $\varphi$ is lower semi-continuous, $\Lambda_\varphi$ is closed.
\item $\pi_\varphi$ is a non-degenerate representation of $B$ on
      $H_\varphi$ such that $\pi_\varphi(a)\,\Lambda_\varphi(b) = \Lambda_\varphi(ab)$
      for every $a \in M(B)$ and $b \in {\cal N}_\varphi$.  (The non-degeneracy of  $\pi_\varphi$ is a consequence of the lower semi-continuity of $\varphi$.)
\end{itemize}

\section{A natural extension of left invariant weights} \label{art1}

In this section, we consider a Hilbert space $H$, a non-degenerate sub-\cst-algebra $B$ of
$B(H)$ and a non-degenerate $^*$-homomorphism $\de$ from $B$ into $M(B \ot B)$ satisfying
the following properties:
\begin{enumerate}
\item $(\de \ot \io)\de = (\io \ot \de)\de$
\item Let $b$ be an element in $M(B)$. If $\de(b) = b \ot 1$, then $b$ belongs to $\C \, 1$.
\end{enumerate}

Furthermore, we assume that $(B,\de)$ arises from a unitary
$W$ on $H \ot H$, i.e. $W$ is a unitary operator on $H \ot H$
such that
\begin{enumerate}
\item $\de(x) = W^* (1 \ot x) W$ for every $x \in B$.
\item $B$ is the closure of the set \  $\{\,(\io \ot \om)(W) \mid \om \in B_0(H)^* \,\}$  \ in $B(H)$.
\end{enumerate}

Consider $v,w \in H$. Then $\om_{v,w}$ will denote the element in $B^*$ such that
$\om_{v,w}(x) = \langle x v , w \rangle$ for every $x \in B$.

\medskip

We will  moreover assume the existence  of a densely defined lower semi-continuous weight
$\vfi$ on $B$ such that we have for every $x \in \Mfi^+$ and every $v \in H$ that
$(\om_{v,v} \ot \io)\de(x)$ belongs to $\Mfi^+$ and $$\vfi\bigl( (\om_{v,v} \ot
\io)\de(x)\bigr) = \vfi(x) \, \langle v , v \rangle \ .$$

\bigskip

The purpose of this section is to arrive at a natural extension of this weight $\vfi$. We
use this extension to prove some interesting results about \cst-algebraic quantum groups
in the last section.

\bigskip

Notice that the left invariance property implies also that $(\om_{v,w} \ot \io)\de(x)$
belongs to $B$ for every $x \in B$ and every $v,w \in H$.

\bigskip

The proof of the following lemma  is essentially the same as the proof of lemma 1.4 of
\cite{Haa}. We include it for the sake of completeness.

\begin{lemma}
Consider $b \in B^+$ such that $(\om_{v,v} \ot \io)\de(b)$ belongs to $\Mfi^+$ for every
$v \in H$. Then there exists a unique element $T$ in $B(H)^+$ such that $\langle T v , v
\rangle  = \vfi\bigl( (\om_{v,v} \ot \io)\de(b) \bigr)$ for every $v \in H$.
\end{lemma}
\begin{demo}
By polarisation, we have for every $v,w \in H$ that $(\om_{v,w} \ot \io)\de(b)$ belongs to
$\Mfi$. This allows us to define a semi-innerproduct $(\, , \,)$ on $H$ such that $(v,w) =
\vfi\bigl((\om_{v,w} \ot \io)\de(b)\bigr)$ for every $v,w \in H$.

Define the subspace $N = \{\, v \in H \mid (v,v) = 0 \,\}$ and turn $\frac{H}{N}$ into a
innerproduct space in the usual way. Define $K$ to be the completion of $\frac{H}{N}$.

Next, we denote the quotient mapping from $H$ into $\frac{H}{N}$ by $F$, so $F$ is a
linear mapping from $H$ into $K$ such that $\langle F(v) , F(w) \rangle = \vfi\bigl(
(\om_{v,w} \ot \io)\de(b)\bigr)$ for every $v,w \in H$. In the next part, we will prove
that $F$ is continuous.

\medskip

Choose a sequence $(v_n)_{n=1}^\infty$ in $H$, $v \in H$ and $w \in K$ such that $(v_n)_{n=1}^\infty \rightarrow v$ and
$(F(v_n))_{n=1}^\infty \rightarrow w$.

Take $\vep > 0$. Then there exists $n_0 \in \N$ such that
$\|F(v_n) - F(v_m)\| \leq \vep$ for every $m,n \in \N$ with $m,n \geq n_0$. This implies that
$$ \vfi\bigl( (\om_{v_n-v_m,v_n-v_m} \ot \io)\de(b) \bigr)
=  \|F(v_n - v_m)\|^2 \leq \vep^2  \hspace{2cm} \text{(*)}$$
for every $m,n \in \N$ with $m,n \geq n_0$.

Take $l \in \N$ with $l \geq n_0$. Because $\bigl(\,(\om_{v_l-v_m,v_l-v_m} \ot
\io)\de(b)\,\bigr)_{m=1}^\infty$ converges to $(\om_{v_l-v,v_l-v} \ot \io)\de(b)$, the lower
semi-continuity of $\vfi$ and the previous inequality imply that
$$\vfi\bigl((\om_{v_l-v,v_l-v} \ot \io)\de(b)\bigr) \leq \vep^2 \, .$$
Consequently, $\|F(v_l) - F(v)\| \leq \vep$.

So we have proven that $(F(v_n))_{n=1}^\infty$ converges to $F(v)$, which implies that
$F(v)=w$.

\medskip

Hence, we see that $F$ is closed implying that $F$ is continuous by the closed graph theorem.

Define $T = F^* F$, then $T$ belongs to $B(H)^+$ and $\langle T v , v \rangle = \langle
F(v) , F(v) \rangle = \vfi\bigl( (\om_{v,v} \ot \io)\de(b)\bigr)$ for every $v \in H$.
\end{demo}

\medskip

We will use the technique of the previous lemma once more in the proof of the next lemma.

\begin{lemma} \label{lem2}
Consider $b \in B^+$ such that $(\om_{v,v} \ot \io)\de(b)$ belongs to $\Mfi^+$ for every
$v \in H$. Let $T$ be the unique element in $B(H)^+$ such that $\langle T v , v \rangle
=$  $\vfi\bigl( (\om_{v,v} \ot \io)\de(b) \bigr)$ for every $v \in H$. Then we have for every $u \in
H \ot H$ that $(\om_{u,u} \ot \io)(\de(b)_{23})$ belongs to $\Mfi^+$ and
$\vfi\bigl((\om_{u,u} \ot \io)(\de(b)_{23})\bigr) = \langle (1 \ot T) \, u , u \rangle$.
\end{lemma}
\begin{demo}
We have by polarisation for every $v,w \in H$ that $(\om_{v,w} \ot
\io)\de(b)$ belongs to $\Mfi$ and
$$\vfi\bigl((\om_{v,w} \ot \io)\de(b)\bigr) = \langle T v , w \rangle \ .$$

Using this, it is not so difficult to check  for every $x \in H \od H$ that $(\om_{x,x}
\ot \io)(\de(b)_{23})$ belongs to $\Mfi^+$ and
$$\vfi\bigl((\om_{x,x} \ot \io)(\de(b)_{23})\bigr) = \langle (1 \ot T) \, x , x \rangle \hspace{2cm} (*)$$

\begin{list}{}{\setlength{\leftmargin}{.4 cm}}

\item Choose $y \in H \ot H$. Then there exists a sequence $(y_n)_{n=1}^\infty$ in $H \od H$ such that $(y_n)_{n=1}^\infty$ converges to $y$.
Then the sequence $\bigl(\,(\om_{y_n,y_n} \ot \io)(\de(b)_{23})\,\bigr)_{n=1}^\infty$
converges to $(\om_{y,y} \ot \io)(\de(b)_{23})$ so the lower semi-continuity of $\vfi$
implies that
\begin{eqnarray*}
\vfi\bigl((\om_{y,y} \ot \io)(\de(b)_{23})\bigr)
& \leq & \liminf \, \bigl(\,\vfi(\,(\om_{y_n,y_n} \ot
\io)(\de(b)_{23})\,)\,\bigr)_{n=1}^\infty \\
& = & \liminf \, \bigl(\,\langle (1 \ot T) \, y_n , y_n \rangle\,\bigr)_{n=1}^\infty
= \langle (1 \ot T) \, y , y \rangle  \ ,
\end{eqnarray*}
where we used equality (*) in the second last equality. So we see that
 $(\om_{y,y} \ot \io)(\de(b)_{23})$ belongs to $\Mfi^+$.

\end{list}

By the discussion above, we know already that $(\om_{u,u} \ot \io)(\de(b)_{23})$ belongs
to $\Mfi^+$.

Take a sequence $(u_n)_{n=1}^\infty$ in $H \od H$ such that $(u_n)_{n=1}^\infty$ converges
to $u$.

\medskip

Choose $\vep > 0$. Then there exists $n_0 \in \N$ such that
$\langle (1 \ot T) (u_n - u_m) , u_n - u_m \rangle \leq \vep^2$ for every
$m,n \in \N$ with $m,n \geq n_0$.

Take $l \in \N$ with $l \geq n_0$. By (*), we have for every $m \in \N$ with $m \geq n_0$
that
$$\vfi\bigl((\om_{u_l-u_m,u_l-u_m} \ot \io)(\de(b)_{23})\bigr)
= \langle (1 \ot T)(u_l-u_m) , u_l-u_m \rangle \leq \vep^2 \ . $$

Because $\bigl(\,(\om_{u_l-u_m,u_l-u_m} \ot \io)(\de(b)_{23})\,\bigr)_{m=1}^\infty$
converges to $(\om_{u_l-u,u_l-u} \ot \io)(\de(b)_{23})$, the lower semi-continuity of
$\vfi$ and the previous inequality imply that
$$\vfi\bigl((\om_{u_l-u,u_l-u} \ot \io)(\de(b)_{23})\bigr) \leq \vep^2 \ .$$

\medskip

Using the fact that the mapping $H \ot H \rightarrow \R^+ : y \mapsto
\vfi\bigl((\om_{y,y} \ot \io)(\de(b)_{23})\bigr)^\frac{1}{2}$ is a semi-norm on
$H \ot H$, this last inequality implies that
\begin{eqnarray*}
& & | \, \vfi\bigl((\om_{u_l,u_l} \ot \io)(\de(b)_{23})\bigr)^\frac{1}{2}
- \vfi\bigl((\om_{u,u} \ot \io)(\de(b)_{23})\bigr)^\frac{1}{2} \, | \\
& & \spat \leq \vfi\bigl((\om_{u_l-u,u_l-u} \ot \io)(\de(b)_{23})\bigr)^\frac{1}{2} \leq
\vep \ .
\end{eqnarray*}

Therefore, we get that the sequence $\bigl(\,\vfi\bigl((\om_{u_n,u_n} \ot
\io)(\de(b)_{23})\bigr)\,\bigr)_{n=1}^\infty$ converges to $\vfi\bigl((\om_{u,u} \ot \io)(\de(b)_{23})\bigr)$.

By (*), we have for every $n \in
\N$ that $$\vfi\bigl((\om_{u_n,u_n} \ot \io)(\de(b)_{23})\bigr) = \langle (1 \ot T)\, u_n , u_n
\rangle$$
which implies  that $\bigl(\,\vfi\bigl((\om_{u_n,u_n} \ot
\io)(\de(b)_{23})\bigr)\,\bigr)_{n=1}^\infty$ converges to $\langle (1 \ot T)\,u , u \rangle$.
Comparing these two results, we see that
$$\vfi\bigl((\om_{u,u} \ot \io)(\de(b)_{23})\bigr) = \langle (1 \ot T) \, u , u \rangle \ .$$
\end{demo}

\medskip

\begin{lemma}  \label{lem1}
Consider $b \in B^+$ such that $(\om_{v,v} \ot \io)\de(b)$ belongs to $\Mfi^+$ for every
$v \in H$. Then there exists a unique positive number $r$ such that  $\vfi\bigl((\om_{v,v}
\ot \io)\de(b)\bigr) = r \, \langle v , v \rangle$ for every $v \in H$.
\end{lemma}
\begin{demo}
Call $T$ the unique element in $B(H)^+$ such that $\vfi\bigl((\om_{v,v} \ot
\io)\de(b)\bigr) = \langle T v , v \rangle$ for every $v \in H$.

Take $v_1,v_2,w_1,w_2 \in H$. By the previous lemma, we know that the element
\newline $(\om_{W(v_1 \ot w_1),W(v_2 \ot w_2)} \ot \io)(\de(b)_{23})$ belongs to $\Mfi$ and
\begin{eqnarray*}
& & \vfi\bigl((\om_{W(v_1 \ot w_1),W(v_2 \ot w_2)} \ot \io)(\de(b)_{23})\bigr)
\\ & & \spat = \langle (1 \ot T) W (v_1 \ot w_1) , W (v_2 \ot w_2) \rangle \\ & & \spat = \langle W^* (1
\ot T) W (v_1 \ot w_1) , v_2 \ot w_2 \rangle \hspace{1.2cm} \text{(a)}
\end{eqnarray*}

Using the fact that $\de(x) = W^* (1 \ot x) W$ for every $x \in B$, we get that
\begin{eqnarray*}
& & (\om_{W(v_1 \ot w_1),W(v_2 \ot w_2)} \ot \io)(\de(b)_{23})\\ & & \spat =  (\om_{v_1
\ot w_1,v_2 \ot w_2} \ot \io)(W^*_{12} \de(b)_{23} W_{12}) \\ & & \spat = (\om_{v_1 \ot
w_1,v_2 \ot w_2} \ot \io)(\,(\de \ot \io)\de(b) \,) \\ & & \spat = (\om_{v_1,v_2} \ot
\om_{w_1,w_2} \ot \io)(\,(\io \ot \de)\de(b) \,) \\ & & \spat = (\om_{w_1,w_2} \ot
\io)\de((\om_{v_1,v_2} \ot \io)\de(b))
\hspace{1.5cm} \text{(b)}
\end{eqnarray*}
By assumption, we have that $(\om_{v_1,v_2} \ot \io)\de(b)$ belongs to $\Mfi$ so the left invariance of $\vfi$ implies that
$(\om_{w_1,w_2} \ot \io)\de((\om_{v_1,v_2} \ot \io)\de(b))$ belongs to $\Mfi$ and
$$\vfi\bigl((\om_{w_1,w_2} \ot \io)\de((\om_{v_1,v_2} \ot \io)\de(b))\bigr)
= \langle w_1 , w_2 \rangle \,\, \vfi\bigl((\om_{v_1,v_2} \ot \io)\de(b)\bigr) \ . $$
Using (b), this implies that
\begin{eqnarray*}
& & \vfi\bigl((\om_{W(v_1 \ot w_1),W(v_2 \ot w_2)} \ot \io)(\de(b)_{23})\bigr)
= \langle w_1 , w_2 \rangle \, \vfi\bigl((\om_{v_1,v_2} \ot \io)\de(b)\bigr) \\
& & \spat = \langle w_1 , w_2 \rangle \, \langle T v_1 , v_2 \rangle
= \langle (T \ot 1) (v_1 \ot w_1) , v_2 \ot w_2 \rangle
\end{eqnarray*}
Comparing this equality with equality (a), we see that
$$\langle W^* (1 \ot T) W (v_1 \ot w_1) , v_2 \ot w_2 \rangle = \langle (T \ot 1) (v_1 \ot w_1) , v_2 \ot w_2 \rangle \ . $$

So we arrive at the conclusion that $W^* (1 \ot T) W  = T \ot 1$.

\medskip

This implies that $(1 \ot T) W = W (T \ot 1)$. So we have for every $\om \in B_0(H)^*$
that $ (\io \ot \om)(W) \,\, T = (\io \ot \om T)(W)$. Because we assumed that $B$ is the
closure of the set \ $\{\,(\io \ot \om)(W) \mid \om \in B_0(H)^* \,\}$ \ in $B(H)$, we
arrive at the conclusion that $B T \subseteq B$. The selfadjointness of $T$ implies that
also $B \, T \subseteq B$, so we see that $T$ is an element of $M(B)$.

Furthermore, the equation $W^* (1 \ot T) W  = T \ot 1$ implies that $\de(T) = T \ot 1$. By assumption, we get the existence of an element $r \in \R^+$ such that $T = r \, 1$.
The lemma follows.
\end{demo}

This lemma justifies the following definition.

\begin{definition} \label{def1}
We define the mapping $\th$ from $B^+$ into $[0,\infty]$ such that we have for every $b \in B^+$ that
\begin{itemize}
\item If $(\om_{v,v} \ot \io)\de(b)$ belongs to $\Mfi^+$ for every $v \in H$, we define $\th(b)$ to be the unique positive number such that
$\vfi\bigl((\om_{v,v} \ot \io)\de(b)\bigr) = \th(b) \, \langle v , v \rangle$ for every $v
\in H$.
\item If there exists $v \in H$ such that $(\om_{v,v} \ot \io)\de(b)$ does not belong to $\Mfi^+$, we define $\th(b) = \infty$.
\end{itemize}
\end{definition}

\begin{proposition}
The mapping $\th$ is a densely defined lower semi-continuous weight on $B$ which extends $\vfi$.
\end{proposition}
\begin{demo}
It is easy to check that $\th$ is a weight on $B$ which extends $\vfi$.
We turn to the lower semi-continuity.

Choose $\al \in \R^+$. Take a sequence $(b_n)_{n=1}^\infty$ in $B^+$ and $b \in B^+$ such
that $(b_n)_{n=1}^\infty \rightarrow b$ and $\th(b_n) \leq \al$ for every $n \in \N$.

Take $w \in H$. Choose $m \in \N$. Because $\th(b_m) \leq \al$, we have by definition that $(\om_{w,w} \ot \io)\de(b_m)$ belongs to $\Mfi^+$
and
$$\vfi\bigl((\om_{w,w} \ot \io)\de(b_m)\bigr) = \th(b_m) \, \langle w , w \rangle
\leq \al \, \langle w , w \rangle \ \ . $$
Because $\bigl(\,(\om_{w,w} \ot \io)\de(b_n)\,\bigr)_{n=1}^\infty$ converges to
$(\om_{w,w} \ot \io)\de(b)$ and $\vfi$ is lower semi-continuous, this implies that
$\vfi\bigl((\om_{w,w} \ot \io)\de(b)\bigr) \leq \al \, \langle w , w \rangle$. In
particular, we get that $(\om_{w,w} \ot \io)\de(b)$ belongs to $\Mfi^+$.

\medskip

Using the definition of $\th$, we see that
$$\th(b) \, \langle v , v \rangle
= \vfi\bigl((\om_{v,v} \ot \io)\de(b)\bigr) \leq \al \, \langle v , v \rangle $$
for every $v \in H$. Consequently, $\th(b) \leq \al$.

\medskip

From this all, the lower semi-continuity of $\th$ follows.
\end{demo}

Because $\th$ is an extension of $\vfi$, we get immediately the following result.

\begin{proposition}
Consider $b \in \cM_\th^+$ and $v \in H$. Then $(\om_{v,v} \ot \io)\de(b)$ belongs to
$\cM_\th^+$ and $$\th\bigl((\om_{v,v} \ot \io)\de(b)\bigr) = \th(b) \, \langle v , v
\rangle \ . $$
\end{proposition}

We can even do better.

\begin{proposition}
Consider $b \in B^+$ such that $(\om_{v,v} \ot \io)\de(b)$ belongs to $\cM_\th^+$ for
every $v \in H$. Then $b$ belongs to $\cM_\th^+$.
\end{proposition}
\begin{demo}
For every $u_1,u_2 \in H \ot H$, we define the element $f_{u_1,u_2} \in B ^*$ such that
$f_{u_1,u_2}(x) = \langle \de(x) \, u_1 , u_2 \rangle$ for every $x \in B$. By lemma
\ref{lem1} (applied to $\th$ instead of $\vfi$), we have the existence of a positive
number $r$ such that $\th\bigl((\om_{v,w} \ot \io)\de(b)\bigr) = r \, \langle v , w
\rangle$ for every $v,w \in H$.

\medskip

Take $v_1,v_2,w_1,w_2 \in H$. By supposition, we have that $(\om_{v_1,w_1} \ot
\io)\de(b)$ belongs to $\cM_\th$. This implies by the definition of $\th$ that
$(\om_{v_2,w_2} \ot \io)\de((\om_{v_1,w_1} \ot \io)\de(b))$ belongs to $\Mfi$ and
\begin{eqnarray*}
& & \vfi\bigl((\om_{v_2,w_2} \ot \io)\de((\om_{v_1,w_1} \ot \io)\de(b))\bigr)
=  \th\bigl((\om_{v_1,w_1} \ot \io)\de(b)\bigr) \,\, \langle v_2 , w_2 \rangle  \\
& & \spat = r \, \langle v_1 , w_1 \rangle \, \langle  v_2 , w_2 \rangle
=  r \, \langle v_1 \ot v_2 , w_1 \ot w_2 \rangle \ .
\end{eqnarray*}
Using the coassociativity of $\de$, we have also that
\begin{eqnarray*}
(\om_{v_2,w_2} \ot \io)\de((\om_{v_1,w_1} \ot \io)\de(b)) & = &
( (\om_{v_1,w_1} \ot \om_{v_2,w_2})\de \ot \io)\de(b) \\
& = & (f_{v_1 \ot w_1,v_2 \ot w_2} \ot \io)\de(b)  \ .
\end{eqnarray*}
So we see that $(f_{v_1 \ot w_1,v_2 \ot w_2} \ot \io)\de(b)$
belongs to $\Mfi$ and that
$$\vfi\bigl((f_{v_1 \ot w_1,v_2 \ot w_2} \ot \io)\de(b)\bigr)
= r \, \langle v_1 \ot w_1 , v_2 \ot w_2 \rangle$$

This implies for every $u \in H \od H$ that $(f_{u,u} \ot \io)\de(b)$ belongs to $\Mfi^+$
and $\vfi\bigl((f_{u,u} \ot \io)\de(b)\bigr) = r \, \langle u , u \rangle$.

This allows, just as in the proof of lemma \ref{lem2}, to conclude for every $u \in H \ot
H$ that $(f_{u,u} \ot \io)\de(b)$ belongs to $\Mfi^+$ and $\vfi\bigl((f_{u,u} \ot
\io)\de(b)\bigr) = r \, \langle u , u \rangle$.

\medskip

Let us now choose $v \in H$. Take $w \in H$ such that $\langle w , w \rangle = 1$. By the
preceding result, we know that $(f_{W^* (w \ot v), W^* (w \ot v)} \ot \io)\de(b)$ belongs
to $\Mfi^+$. Using the fact that $\de(x) = W^* (1 \ot x) W$ for every $x \in B$, it is
easy to check that $f_{W^* (w \ot v), W^* (w \ot v)} = \om_{v,v}$ so we get that
$(\om_{v,v} \ot \io)\de(b)$ belongs to $\Mfi^+$.

By definition, we get that $b$ belongs to $\cM_\th^+$.
\end{demo}

\section{A useful result concerning KMS-weights}

We will use the following definition of a KMS-weight (see \cite{JK1}):

\begin{definition}
Consider a \cst-algebra $B$ and a densely defined lower semi-continuous weight $\vfi$ on
$B$ such that there exist a norm-continuous one parameter group $\si$ on $B$ such that
\begin{enumerate}
\item We have that $\vfi \si_t = \vfi$ for every $t \in \R$.
\item We have that $\vfi(x^* x) = \vfi(\si_\frac{i}{2}(x)  \si_\frac{i}{2}(x)^*)$
for every  $x \in D(\si_\frac{i}{2})$.
\end{enumerate}
Then $\vfi$ is called a KMS-weight on $B$ and $\si$ is called a modular group for $\vfi$.
\end{definition}

This is not the usual definition of a KMS-weight (see \cite{Comb})  but we prove in
\cite{JK1} that this definition is equivalent to the usual one.

\medskip

A rather useful proposition concerning KMS-weights is the followoing one (see \cite{JK1})
:

\begin{proposition}
Consider a  \cst-algebra $B$ and let $\vfi$ be a KMS-weight on $B$ with modular group
$\si$. Consider moreover a GNS-construction $(H_\vfi,\lafi,\pifi)$ for $\vfi$. Then :
\begin{itemize}
\item There exists a unique anti-unitary operator $J$ on $H$ such that
$J \lafi(x) = \lafi(\si_{\frac{i}{2}}(x)^*)$ for every $x \in \Nfi \cap
D(\si_{\frac{i}{2}})$
\item  We have for every $x \in \Nfi$ and every $a \in D(\si_{\frac{i}{2}})$
that $x a$ belongs to $\Nfi$ and $$\lafi(x a) = J \pifi(\si_{\frac{i}{2}}(a))^* J
\lafi(x) \ .$$
\end{itemize}
\end{proposition}

\medskip

We now prove that a KMS-weight has no proper extensions which are relatively invariant
under its modular group.

\begin{proposition} \label{prop2}
Consider a KMS-weight $\vfi$ on a \cst-algebra $B$ with modular group $\si$. Let $\eta$ be
a lower semi-continuous weight on $B$ which is an extension of $\vfi$ and such that $\eta$
is relatively invariant under $\si$. Then $\vfi = \eta$.
\end{proposition}
\begin{demo}
Take a GNS-construction $(H_\vfi,\lafi,\pifi)$ for $\vfi$ and a
GNS-construction $(H_\eta,\la_\eta,\pi_\eta)$ for $\eta$.

Because $\eta$ is relatively invariant under $\si$, we get the existence of a strictly
positive number $\lambda$ such that $\eta \si_t = \lambda^t \, \eta$ for every $t \in \R$.

So we get the existence a positive injective operator $T$ in $H_\eta$ such that $T^{it}
\la_\eta(a) = \lambda^{-\frac{t}{2}}  \, \la_\eta(\si_t(a))$ for every $a \in \cN_\eta$ and $t \in \R$.

\medskip

Choose $y \in \cN_\eta$.

Define for every $n \in \N$ the element
$$y_n = \frac{n}{\sqrt{\pi}} \int \exp(-n^2 t^2) \, \si_t(y) \, dt $$
which is clearly analytic with respect to $\si$. We have also that
$(y_n)_{n=1}^\infty$ converges to $y$.

\medskip

By \cite{JK1}, we have for every $n \in \N$ that $y_n$ belongs to $\cN_\eta$ and
$$\la_\eta(y_n) = \frac{n}{\sqrt{\pi}} \int \exp(-n^2 t^2) \,
\lambda^\frac{t}{2} \, T^{it} \la(y) \, dt  \ .$$
This implies immediately that $(\la_\eta(y_n))_{n=1}^\infty$ converges to $\la_\eta(y)$.

\medskip

We can also take an approximate unit $(e_i)_{i \in I}$ for $B$ in $\Nfi$.
Then  $(e_i \, y_n)_{(i,n) \in I \times \N}$ converges to $y$.

We have also for every $i \in I$ and $n \in \N$ that $e_i \, y_n$ belongs to $\cN_\eta$
and $\la_\eta(e_i \, y_n) = \pi_\eta(e_i) \la_\eta(y_n)$. Consequently, the net
$(\la_\eta(e_i \, y_n))_{(i,n) \in I \times \N}$ converges to $\la_\eta(y)$.

\medskip

We have for every $i \in I$ and $n \in \N$ that $e_i$ belongs to $\Nfi$ and $y_n$ belongs
to $D(\si_\frac{i}{2})$ implying that $e_i \, y_n$ belongs to $\Nfi$ by the previous
proposition.

\medskip

Because $\vfi \subseteq \eta$, we have moreover for every $i,j \in I$ and $m,n \in \N$ that
$$\|\lafi(e_i \, y_n) - \lafi(e_j \, y_m)\| = \|\la_\eta(e_i \, y_n) - \la_\eta(e_j \, y_m)\|  \ .$$
This last equality  implies that the net $(\lafi(e_i \, y_n))_{(i,n) \in I \times \N}$ is Cauchy and hence convergent in $H_\vfi$. Therefore, the closedness of $\lafi$ implies that $y$ is an element of $\cN_\vfi$.
The proposition follows.
\end{demo}

\section{Improving the invariance of left Haar weights}

In \cite{Kus}, we constructed the reduced \cst-algebraic quantum group out of an algebraic
one. We constructed a left invariant weight on this \cst-algebraic quantum group. In this
section, we will improve the invariance property of this left Haar weight. We will also
prove a uniqueness property. It is very well possible that the techniques used in this
section can be useful in the more case general case of \cst-algebraic quantum groups
according to Masuda, Nakagami \& Woronowicz (upcoming paper).

\bigskip

We will start this section with a small overview concerning algebraic quantum groups and
the \cst-algebraic quantum groups arising from them. For a more detailed treatment, we
refer to \cite{Kus}.

\medskip

So consider an algebraic quantumgroup $(A,\de)$ according to A. Van Daele (see
\cite{VD1}). This means that $(A,\de)$ satisfies the following properties :

\medskip

The object $A$ is  a non-degenerate $^*$-algebra  and $\de$ is a  non-degenerate
$^*$-homomorphism \newline $\de : A \rightarrow M(A \od A)$ satisfying the following
properties:
\begin{enumerate}
\item  $(\de \od \io) \de = (\io \od \de)\de$.
\item The linear mappings $T_1$,$T_2$ from $A \od A$ into $M(A \od
A)$ such that
\begin{eqnarray*}
& & T_1(a \ot b) = \de(a)(1 \ot b) \hspace{1.5cm} \text{ and }   \\
& & T_2(a \ot b) = \de(a)(b \ot 1)
\end{eqnarray*}
for all $a,b \in A$, are bijections from $A \od A$ into $A \od A$.
\end{enumerate}

We assume furthermore the existence of a non-trivial positive linear
functional $\vfi$ on $A$ such that we have for all  $a,b \in A$ that
$$ (\io \od \vfi)(\de(a)(b \ot 1))
= \vfi(a) \, b \ . $$

\bigskip

Next, we take a GNS-pair $(H,\la)$ of the left Haar functional $\vfi$ on $A$. This means
that $H$ is a Hilbert space and that $\la$ is a linear mapping from $A$ into $H$ such that
\begin{enumerate}
\item The set $\la(A)$ is dense in $H$.
\item We have for every $a,b \in A$ that $\langle \la(a) , \la(b) \rangle
= \vfi(b^* a)$.
\end{enumerate}

\medskip

As usual, we can associate a multiplicative unitary to $(A,\de)$ :

\begin{definition}
We define $W$ as the unique unitary element in $B(H \ot H)$ such that
\newline $W \, (\la \od \la)(\de(b)(a \ot 1)) = \la(a) \ot \la(b)$ for every $a,b \in A$. The element $W$ is called the fundamental unitary associated to $(A,\de)$.
\end{definition}

\bigskip

The GNS-pair $(H,\la)$ allows us to represent $A$ by bounded operators
on $H$ :

\begin{definition}
We define $\pi$ as the unique $^*$-homomorphism from $A$ into $B(H)$ such that $\pi(a)
\la(b) = \la(a b)$ for every $a,b \in A$. We have also that $\pi$ is injective.
\end{definition}

\medskip

Notice also that it is not immediate that $\pi(x)$ is a bounded operator on $H$ (because
$\vfi$ is merely a functional, not a weight), but the boundedness of $\pi(x)$ is connected
with the following equality :

\medskip

We have for every $a,b \in A$ that
\begin{equation}
\pi\bigl((\io \ot \vfi)( \de(b^*)(1 \ot a))\bigr)
= (\io \ot \om_{\la(a),\la(b)})(W)      \label{eq1.1}
\end{equation}

The mapping $\pi$ makes it possible to define our reduced
\cst-algebra :

\begin{definition}
We define $\ar$ as the closure of $\pi(A)$ in $\op$. So $\ar$ is a non-degenerate
sub-\cst-algebra of $B(H)$.
\end{definition}

Equation \ref{eq1.1} implies that
$$
\ar = \text{\ \ closure of \ \ }
\{\, (\io \ot \om)(W) \mid \om \in \cop^* \,\}  \text{\ \ in \  }  \op \ .
$$

As usual, we use the fundamental unitary to define a comultiplication on $A_r$. We will
denote it by $\de_r$.

\begin{definition}
We define the mapping $\de_r$ from $\ar$ into $B(H \ot H)$ such that $\de_r(x) =  W^* (1
\ot x) W$ for all $x \in \ar$. Then $\de_r$ is an injective $^*$-homomorphism.
\end{definition}

It is not so difficult to show that $\de_r$ on the \cst-algebra level is an extension of
$\de$ on the $^*$-algebra level :

\begin{result}
We have for all $a \in A$ and $x \in A \od A$ that $(\pi \od \pi)(x) \, \de_r(\pi(a)) =
(\pi \od \pi)(x \, \de(a))$ and $\de_r(\pi(a)) \, (\pi \od \pi)(x)
= (\pi \od \pi)(\de(a)\,x).$
\end{result}

\medskip

Using the above result, it is not so hard to prove the following theorem :

\begin{theorem}
We have that $A_r$ is a non-degenerate sub-\cst-algebra of $\op$ and $\de_r$ is a
non-degenerate injective $^*$-homomorphism from $A_r$ into $M(A_r \ot A_r)$ such that :
\begin{enumerate}
\item $(\de_r \ot \io)\de_r = (\io \ot \de_r)\de_r$
\item The vector spaces $\de_r(A_r)(A_r \ot 1)$ and $\de_r(A_r)(1 \ot A_r)$ are dense subspaces
      of $A_r \ot A_r$.
\end{enumerate}
\end{theorem}

\bigskip

By lemma 7.11 of \cite{Kus}, we have the following result (which is in fact connected with
the existence of the left invariant weight.)

\begin{result}
Consider an element $x \in M(A_r)$. If $\de_r(x) = x \ot 1$, then $x$ belongs to $\C \, 1$.
\end{result}

\medskip

So we see that the pair $(A_r,\de_r)$ satisfies the requirements of the beginning of
section \ref{art1}.

\bigskip

In \cite{Kus}, we use the functional $\vfi$ to construct a left invariant weight on $A_r$.
This weight will be denoted by the symbol $\vfi_r$ (in \cite{Kus}, it was denoted by
$\vfi$).

\medskip

The  weight is determined by the following theorem (see \cite{Kus},  theorem 6.12,
proposition 6.2 and the remarks after it) :

\begin{theorem}      \label{thm3.1}
There exists a unique closed linear map $\lar$ from  within $A_r$ into $H$ such that
$\pi(A)$ is a core for $\lar$ and $\lar(\pi(a)) = \la(a)$ for every $a \in A$.

There exists moreover a unique weight $\vfi_r$ on $A_r$ such that
$(H,\lar,\io)$ is a GNS-construction for $\vfi_r$.

We have also that $\pi(A) \subseteq \Mfir$ and that $\vfi_r(\pi(a)) =
\vfi(a)$ for every $a \in A$.
\end{theorem}

\medskip

\begin{proposition}
The weight $\vfi_r$ is a faithful KMS-weight. We denote the modular group of $\vfi_r$ by
$\si$.
\end{proposition}

\medskip

This weight satisfies the following left invariance property (corollary 6.14 of
\cite{Kus}).

\begin{proposition}
Consider $a \in \Mfir$ and $\om \in A_r^*$. Then $(\om \ot \io)\de_r(a)$ belongs to $\Mfir$
and $$\vfi_r\bigl((\om \ot \io)\de_r(a)\bigr) = \om(1) \, \vfi_r(a) \ .$$
\end{proposition}

We have even a stronger form of left invariance condition, but we only need this one.

\bigskip

So we see that the weight $\vfi_r$ satisfies the requirements of the beginning of section
\ref{art1}.

\medskip

Now we have enough information to prove some results.

\bigskip

The first theorem improves the left invariance property of $\vfi_r$.

\begin{theorem}
Consider $a \in A_r^+$ such that $(\om_{v,v} \ot \io)\de_r(a)$ belongs to
$\Mfir^+$ for every $v \in H$. Then $a$ belongs to $\Mfir^+$.
\end{theorem}
\begin{demo}
It is clear that we can apply the results from section \ref{art1} to the  weight $\vfi_r$.
So, call $\th$ the weight on $A_r$ arising from $\vfi_r$ as described in definition
\ref{def1}. Hence, $\th$ is a densely defined lower semi-continuous weight on $A_r$ which
extends $\vfi_r$.

\medskip

Choose $b \in \cM_\th^+$ and $t \in \R$. Take $v \in H$.

In definition 5.1 of \cite{Kus}, we defined a $^*$-automorphism $\tau_t$ on $A_r$ which
was implemented by a positive injective (generally unbounded) operator $M$, i.e.
$\tau_t(x) = M^{it} x M^{-it}$ for every $x \in A_r$. In proposition 5.7 of \cite{Kus}, we
proved the commutation $\de_r\, \si_t = (\tau_t \ot \si_t)\de_r$. Therefore,
\begin{eqnarray*}
(\om_{v,v} \ot \io)\de_r(\si_t(b)) & = & (\om_{v,v} \ot \io)((\tau_t \ot
\si_t)\de_r(b)) \\
& = & \si_t((\om_{M^{-it} v , M^{-it} v} \ot \io)\de_r(b)) \ .
\hspace{1.5cm} \text{(*)} \end{eqnarray*}
By the definition of $\th$, we have that $(\om_{M^{-it} v , M^{-it} v} \ot \io)\de_r(b)$
belongs to $\Mfir^+$ and
$$\vfi_r\bigl((\om_{M^{-it} v , M^{-it} v} \ot \io)\de_r(b)\bigr)
= \langle M^{-it} v , M^{-it} v \rangle \, \th(b)
= \langle v , v \rangle \, \th(b) \ . $$
This implies that $\si_t((\om_{M^{-it} v , M^{-it} v} \ot \io)\de_r(b))$ belongs to
$\Mfir^+$ and
$$\vfi_r\bigl(\si_t((\om_{M^{-it} v , M^{-it} v} \ot \io)\de_r(b))\bigr)
= \vfi_r\bigl((\om_{M^{-it} v , M^{-it} v} \ot \io)\de_r(b)\bigr)
 = \langle v , v \rangle \, \th(b) \ . $$
Therefore, equation (*) implies that $(\om_{v,v} \ot \io)\de_r(\si_t(b))$
belongs to $\Mfir^+$ and
$$\vfi_r\bigl((\om_{v,v} \ot \io)\de_r(\si_t(b))\bigr) = \langle v , v \rangle \, \th(b) \ . $$

We get by the definition of $\th$ that $\si_t(b)$ belongs to $\cM_\th^+$
and that
$$\langle v , v \rangle \, \th(\si_t(b))
= \vfi_r\bigl((\om_{v,v} \ot \io)\de_r(\si_t(b))\bigr) = \langle v , v \rangle \, \th(b)$$ for every $v \in H$. So we see that $\th(\si_t(b)) = \th(b)$.

\medskip

Consequently, we have proven that $\th$ is invariant with respect to $\si$. Combining this
with the fact that $\th$ is an extension of $\vfi_r$ and using proposition \ref{prop2}, we
see that $\vfi_r = \th$. The theorem follows.
\end{demo}

\bigskip

Next, we prove some results about the uniqueness of the left Haar weight on the
\cst-algebra level. In both cases, the uniqueness follows from some left invariance
property and some relative invariance property.

\bigskip

We first start with a weaker form of uniqueness. The proof strongly resembles the proof of
proposition 7.1 of \cite{Kus}.

\begin{remark} \rm
If we look at the proof of proposition 6.9 of \cite{Kus}, we see that also the following
weaker form is true:

Consider a dense left ideal $N$ in $A_r$ such that $(\om_{v,v} \ot \io)\de_r(a)$ belongs to $N$ for all $a \in N$ and all $v \in H$.
Then $\pi(A)$ is a subset of $N$.
\end{remark}

\medskip

\begin{proposition} \label{prop1}
Consider  a densely defined lower semi-continuous weight $\eta$ on $A_r$ such that we have
for every $a \in \cM_\eta^+$ and every $v \in H$ that $(\om_{v,v} \ot \io)\de_r(a)$ belongs
to $\cM_\eta^+$ and $$\eta\bigl((\om_{v,v} \ot \io)\de_r(a)\bigr) = \langle v , v \rangle \,
\eta(a) \ .$$
Then there exists a positive number $r$ such that $r \, \vfi_r \subseteq \eta$.
\end{proposition}
\begin{demo}
We know that $\cN_\eta$ is a dense left ideal in $A$.
Choose $a \in \cN_\eta$ and $v \in H$, then
$$ [(\om_{v,v} \ot \io)\de_r(a)]^* \, [(\om_{v,v} \ot \io)\de_r(a)]
\leq \|\om_{v,v}\| \,\, (\om_{v,v} \ot \io)\de_r(a^* a) \ . $$
Because $a^* a$ belongs to $\cM_\eta^+$,  we have by assumption that
$(\om_{v,v} \ot \io)\de_r(a^* a)$ belongs to $\cM_\eta^+$. Therefore,
the previous inequality implies that  the element
$[(\om_{v,v} \ot \io)\de_r(a)]^* \, [(\om_{v,v} \ot \io)\de_r(a)]$ belongs
to $\cM_\eta^+$. Hence, $(\om_{v,v} \ot \io)\de_r(a)$ belongs to
$\cN_\eta$.

By the remark before this proposition, we can conclude from this all that $\pi(A)$ is a
subset of $\cN_\eta$. Because $A^* A = A$, we get that  $\pi(A)$ is a subset of
$\cM_\eta$.

\medskip

So we can define the positive linear functional $\phi$ on $A$ such that $\phi(a) =
\eta(\pi(a))$ for every $a \in A$.

Choose $a,b \in A$. Take $c \in B$. Then
\begin{eqnarray*}
\vfi(\,c^* \, (\io\od \phi)(\de(a)(b \ot 1))\,)
& = & \phi\bigl((\vfi \od \io)((c^* \ot 1)\de(a)(b \ot 1)) \bigr) \\
& = & \eta\bigl(\,\pi(\,(\vfi \od \io)((c^* \ot 1)\de(a)(b \ot 1)) \, )\, \bigr) \\
& = & \eta\bigl( (\om_{\la(b),\la(c)} \ot \io)\de_r(\pi(a))\bigr) \ .
\end{eqnarray*}
Because $\pi(a)$ is an element of $\cM_\eta$, we get by assumption
that  $(\om_{\la(b),\la(c)} \ot \io)\de_r(\pi(a))$ belongs to
$\cM_\eta$ and $$ \eta\bigl((\om_{\la(b),\la(c)} \ot \io)\de_r(\pi(a))\bigr) =
\eta(\pi(a)) \, \langle \la(b) , \la(c) \rangle = \phi(a) \, \vfi(c^*
b) \ .$$
So we see that $\vfi(\,c^* \, (\io \od \phi)(\de(a)(b \ot 1))\, )  =
\phi(a) \, \vfi(c^* b)$.

Hence, the faithfulness of $\vfi$ (see proposition 3.4 of \cite{VD1}) implies that $(\io
\od
\phi)(\de(a)(b \ot 1)) = \phi(a) \, b$.

\medskip

Therefore, $\phi$ is a left invariant functional on $A$. By the uniqueness of the Haar
functional on the $^*$-algebra level (see theorem 3.7 of \cite{VD1}), we get the existence
of a positive number $r$ such that $\phi = r \, \vfi$.

Consequently, we have proven that $\pi(A) \subseteq \cM_\eta$ and $\eta(x) = r \,
\vfi_r(x)$ for every $x \in \pi(A)$.

\medskip

Take a GNS-construction $(H_\eta,\la_\eta,\pi_\eta)$ for $\eta$.

Choose $y \in \Nfir$. Because $\pi(A)$ is a core for $\lar$ (see theorem \ref{thm3.1}),
there exists a sequence $y_n \in \pi(A)$ such that $(y_n)_{n=1}^\infty$ converges to $y$
and $(\lar(y_n))_{n=1}^\infty$ converges to $\lar(y)$. By the first part of the proof, we
know already for every $n \in \N$ that $y_n$ belongs to $\cN_\eta$. Furthermore, we have
for every $m,n \in \N$ that
\begin{eqnarray*}
& & \|\la_\eta(y_m) - \la_\eta(y_n)\|^2 = \eta( (y_m-y_n)^* (y_m-y_n) )
\\
& & \spat =  r \, \vfi_r( (y_m-y_n)^* (y_m-y_n) )
=  r \, \|\lar(y_m) - \lar(y_n)\|^2 \ .
\end{eqnarray*}
This implies that $(\la_\eta(y_n))_{n=1}^\infty$ is Cauchy and hence convergent. So the
closedness of $\la_\eta$ implies that $y$ belongs to $\cN_\eta$ and that
$(\la_\eta(y_n))_{n=1}^\infty$ converges to $\la_\eta(y)$.

We have for every $n \in \N$ that
$$\langle \la_\eta(y_n) , \la_\eta(y_n) \rangle = \eta(y_n^* y_n)
= r \, \vfi_r(y_n^* y_n) = r \,  \langle \lar(y_n) , \lar(y_n) \rangle $$
which implies that $\langle \la_\eta(y) , \la_\eta(y) \rangle = r \, \langle \lar(y) ,
\lar(y)
\rangle$. The proposition follows.
\end{demo}

\medskip

Under stronger conditions, we can prove a real uniqueness result which will be shown
below.

\medskip

From proposition 3.15 of \cite{Kus}, we know that there exists a unique norm continuous
one-parameter group $K$ on $A_r$ such that $(\si_t \ot K_t)\de_r = \de_r \, \si_t$ for
every $t \in \R$. This one parameter group is used to formulate some uniqueness result.

\begin{theorem}
Consider a non-zero densely defined lower semi-continuous weight $\eta$ on $A_r$ such that
\begin{enumerate}
\item We have for every $a \in \cM_\eta^+$ and $v \in H$ that
$(\om_{v,v} \ot \io)\de_r(a)$ belongs to $\cM_\eta^+$ and
$$\eta\bigl((\om_{v,v} \ot \io)\de_r(a)\bigr) = \eta(a) \, \langle v ,v \rangle \ .$$
\item The weight $\eta$ is relatively invariant under $K$.
\end{enumerate}
Then there exists a unique strictly positive number $r$ such that $\eta = r \, \vfi_r$.
\end{theorem}
\begin{demo}
By assumption, there exists a unique strictly positive number $\lambda$ such that $\eta \,
K_t = \lambda^t \, \eta$ for every $t \in \R$.

The previous proposition implies the existence of a positive number $r$ such that $\eta$
is an extension of $r \, \vfi_r$. We automatically have that $r \neq 0$. (If $r$ would be
0, then $\eta$ would be 0 on the set $\pi(A)$ so the lower semi-continuity of $\eta$ would
imply that $\eta=0$.)

Because $\vfi_r$ is invariant under $\si$, we have also the existence of an injective
positive operator $\nab$ in $H$ such that $\nab^{it} \lar(x) = \lar(\si_t(x))$ for every
$t \in \R$ and $x \in \Nfir$. Then $\si_t(x) = \nab^{it} x \nab^{-it}$ for every $x \in
A_r$ and $t \in \R$.

\medskip

Again, the  requirements of the beginning of  section \ref{art1} are satisfied and we get
a weight $\th$ arising from $\eta$ according to definition \ref{def1}. Hence, $\th$ is a
densely defined lower semi-continuous weight on $A_r$ which extends $\eta$. So $\th$ will
also extend $r \, \vfi_r$.

\medskip

Choose $b \in \cM_\th^+$ and $t \in \R$. Take $v \in H$. By the remarks before this
theorem, we get that
\begin{eqnarray*}
(\om_{v,v} \ot \io)\de_r(\si_t(b)) & = & (\om_{v,v} \ot \io)((\si_t \ot K_t)\de_r(b)) \\
& = & K_t((\om_{\nab^{-it} v , \nab^{-it} v} \ot \io)\de_r(b)) \ .
\hspace{1.5cm} \text{(*)} \end{eqnarray*}
By the definition of $\th$, we have that $(\om_{\nab^{-it} v , \nab^{-it} v} \ot \io)\de_r(b)$ belongs to $\cM_\eta^+$ and
$$\eta\bigl((\om_{\nab^{-it} v , \nab^{-it} v} \ot \io)\de_r(b)\bigr)
= \langle \nab^{-it} v , \nab^{-it} v \rangle \, \th(b)
= \langle v , v \rangle \, \th(b) \ . $$
This implies that $K_t((\om_{\nab^{-it} v , \nab^{-it} v} \ot \io)\de_r(b))$ belongs to $\cM_\eta^+$ and
$$\eta\bigl(K_t((\om_{\nab^{-it} v , \nab^{-it} v} \ot \io)\de_r(b))\bigr)
= \lambda^t \, \eta\bigl((\om_{\nab^{-it} v , \nab^{-it} v} \ot \io)\de_r(b)\bigr) = \lambda^t \, \langle v , v \rangle \, \th(b) \ . $$

Therefore, equation (*) implies that $(\om_{v,v} \ot \io)\de_r(\si_t(b))$
belongs to $\cM_\eta^+$ and
$$\eta\bigl((\om_{v,v} \ot \io)\de_r(\si_t(b))\bigr) = \lambda^t \, \langle v , v \rangle \, \th(b) \ . $$

Therefore we get by the definition of $\th$ that $\si_t(b)$ belongs to $\cM_\th^+$ and that
$$\langle v , v \rangle \, \th(\si_t(b))
= \eta\bigl((\om_{v,v} \ot \io)\de_r(\si_t(b))\bigr) = \lambda^t \, \langle v , v \rangle \, \th(b)$$ for every $v \in H$. So we see that $\th(\si_t(b)) = \lambda^t \, \th(b)$.

\medskip

Consequently, we have proven that $\th$ is relatively invariant with respect to $\si$.
Combining this with the fact that $\th$ is an extension of $r \,
\vfi_r$ and using
proposition \ref{prop2}, we see that $\th = r \, \vfi_r$. Because $r \, \vfi_r \subseteq
\eta \subseteq \th$, we get that $\eta = r \, \vfi_r$.
\end{demo}

Of course, also the following uniqueness result is valid. It follows immediately from
propositions \ref{prop1} and \ref{prop2} (without using the extension $\th$).

\begin{theorem}
Consider a non-zero densely defined lower semi-continuous weight $\eta$ on $A_r$ such that
\begin{enumerate}
\item We have for every $a \in \cM_\eta^+$ and $v \in H$ that
$(\om_{v,v} \ot \io)\de_r(a)$ belongs to $\cM_\eta^+$ and
$$\eta\bigl((\om_{v,v} \ot \io)\de_r(a)\bigr) = \eta(a) \, \langle v ,v \rangle \ .$$
\item The weight $\eta$ is relatively invariant under $\si$.
\end{enumerate}
Then there exists a unique strictly positive number $r$ such that $\eta = r \, \vfi_r$.
\end{theorem}

\end{document}